\documentclass[11pt]{article}
\usepackage{amssymb,amsfonts,amsmath,latexsym,epsf,tikz,url}

\newtheorem{theorem}{Theorem}[section]
\newtheorem{proposition}[theorem]{Proposition}

\newtheorem{corollary}[theorem]{Corollary}
\newtheorem{lemma}[theorem]{Lemma}

\newtheorem{definition}[theorem]{Definition}
\newtheorem{problem}[theorem]{Problem}
\newtheorem{question}[theorem]{Question}

\newcommand{\proof}{\noindent{\bf Proof.\ }}
\newcommand{\qed}{\hfill $\square$\medskip}

\begin{document}

\title{The distinguishing number of groups based on the distinguishing number of subgroups}

\author{
Saeid Alikhani  $^{}$\footnote{Corresponding author}
\and
Samaneh Soltani
}

\date{\today}

\maketitle

\begin{center}
Department of Mathematics, Yazd University, 89195-741, Yazd, Iran\\
{\tt alikhani@yazd.ac.ir, s.soltani1979@gmail.com}
\end{center}


\begin{abstract}
 Let $\Gamma$ be a group acting on a set $X$. The distinguishing number for this action of $\Gamma$ on $X$, denoted by $D_{\Gamma}(X)$, is the smallest natural number $k$ such that the elements of $X$ can be labeled with $k$ labels so that any label-preserving element of $\Gamma$ fixes all $x \in X$. In particular, if the action is faithful, then the only element of $\Gamma$ preserving labels is the identity. In this paper, we  obtain an upper bound on the distinguishing number of a set knowing the distinguishing number of a set under the action of a subgroup. By the concept of motion, we obtain an upper bound for the distinguishing number of a group. Motivated by a  problem 
 (Chan 2006),  we characterize $D_{\Gamma,H}(X)$ which is the smallest number of labels admitting a labeling of $X$ such that the only elements of $\Gamma$ that induce label-preserving permutations lie in $H$. Finally, we state  two algorithms for obtaining an upper and a lower bound for $D_{\Gamma , H}(X)$.
\end{abstract}


\medskip
\noindent{\bf AMS Subj.\ Class.}: 05E15, 05C15, 20B25, 20D60

\section{Introduction }

Let $G=(V,E)$ be a simple connected graph. We use the standard graph notation (\cite{Sandi}). In particular, ${\rm Aut}(G)$ denotes the automorphism group of $G$.  
Distinguishing labeling was first defined by Albertson and Collins \cite{Albert} for graphs. A labeling of a graph G, $\phi : V(G)\rightarrow \{1,2,\ldots ,r\}$, is said to be \textit{$r$-distinguishing} if no nontrivial automorphism of $G$ preserves all the vertex labels. In other words, $\phi$ is $r$-distinguishing if for any $\sigma\in {\rm Aut}(G)$, $\sigma \neq id$, there is a vertex $x$ such that $\phi (x) \neq \phi (\sigma (x))$. The \textit{distinguishing number} of a graph $G$ is defined as 
\begin{equation*}
D(G) ={\rm min}\{r:~ {\rm there~ exists~ an~ r-distinguishing ~labeling ~of~ G}\}.
\end{equation*}
 
The distinguishing number of several families of graphs, including trees, hypercubes, generalized Petersen graphs, friendship and book graphs have been studied in \cite{soltani,Bogstad,Chan,Cheng}, and \cite{Potanka}. The concept was naturally extended to general group actions \cite{Tymockzko}. 
Let $\Gamma$ be a group acting on a set $X$. If $g$ is an element of $\Gamma$ and $x$ is in $X$ then
denote the action of $g$ on $x$ by $g.x$. Write $\Gamma .x$ for the orbit containing $x$. Recall that stabilizer of the subset $Y \subseteq X$ is defined to be ${\rm Stab}_{\Gamma}(Y ) = \{g \in \Gamma : g.y = y ~{\rm for~ all}~ y \in Y \}$. We sometimes omit the subscript and write ${\rm Stab}(Y)$. For a positive integer $r$, an \textit{$r$-labeling} of $X$ is an onto function $\phi : X\rightarrow \{1,2,\ldots ,r\}$. We say $\phi$ is a \textit{distinguishing labeling} (with respect to the action of $\Gamma$) if  the only group elements that preserve the labeling are in ${\rm Stab}_{\Gamma}(X)$. The \textit{distinguishing number} $D_{\Gamma}(X)$ of the action of  $\Gamma$ on $X$ is defined as 
\begin{equation*}
D_{\Gamma}(X) ={\rm min}\{r:~ {\rm there~ exists~ an~ r-distinguishing~ labeling}\}.
\end{equation*}
In particular, if the action is faithful, then the only element of $\Gamma$ preserving labels is the identity. In this paper we suppose that the action of a group on a set is faithful. Tymoczko in \cite{Tymockzko} shows that:
\begin{theorem}{\rm \cite{Tymockzko}}\label{uppdisgroup}
If $|\Gamma|$ is at most $k!$, then $D_{\Gamma}(X)$ is at most $k$.
\end{theorem}
 
For a given group $\Gamma$, Albertson and Collins \cite{Albert} defined the \textit{distinguishing set} of $\Gamma$ graphs as
\begin{equation*}
D(\Gamma) =\{D(G) :~ G ~ {\rm is ~a ~graph~with~ Aut}(G)=\Gamma \}.
\end{equation*}
Wong and Zhu \cite{Wong}  defined the \textit{distinguishing set} of $\Gamma$ actions as
\begin{equation*}
D^*(\Gamma) =\{D_{\Gamma}(X) :~ \Gamma ~ {\rm acts ~faithfully~ on  }~ X \}.
\end{equation*}

Observe that for a graph $G$, $D(G) = D_{{\rm Aut(G)}}(V(G))$, so for any group $\Gamma$, $D(\Gamma)\subseteq D^*(\Gamma)$.  We need the following theorem:

\begin{theorem}{\rm \cite{Albert, Chan2}}\label{abelidihedral}
 If $\Gamma$ is a nontrivial abelian group, then $D(\Gamma)= D^*(\Gamma)=\{2\}$, also if $D_n$ is the dihedral group of order $2n$, then $D(D_n)= D^*(D_n)=\{2\}$ unless $n=3,4,5,6,10$, in which cases $D(D_n)= D^*(D_n)=\{2,3\}$.
\end{theorem}

 A group $\Gamma$ is called \textit{almost simple}  if $T \leqslant \Gamma \leqslant {\rm Aut}(T)$ for some nonabelian simple group $T$. For almost
simple groups, we use the notation of the Atlas \cite{J. H. Conway}. In particular, $L_d(q)$ denotes the projective special linear group of dimension $d$ defined over a $q$-element field, and $M_n$, for $n \in \{11,12,22,23,24\}$, denotes the Mathieau groups.  Seress et. al., have obtained the distinguishing of these groups as follows:

\begin{theorem}{\rm \cite{Seress0}}\label{symmthm}
For all $n\geq 3$,
\begin{itemize}
\item[(i)] $D^*(S_n)=\big\{\lceil  n^{1/k}\rceil: k=1,2, \ldots \big\}\cup \big\{\lceil  (n-1)^{1/k}\rceil : k=1,2, \ldots \big\}$;
\item[(ii)] $D(S_n)=\big\{\lceil n^{1/k}\rceil: k=1,2, \ldots \big\}$;
\item[(iii)] $D^*(A_n)=\big\{\lceil (n-1)^{1/k}\rceil: k=1,2, \ldots \big\}$, except that $D^*(A_5)=\{2,3,4\}$;
\item[(iv)] $D(A_n)=\{2\}$.
\end{itemize}
\end{theorem}
\begin{theorem}{\rm \cite{Seress0}}\label{almostthm}
 If $\Gamma $ is an almost simple group and $\Gamma \neq S_n,A_n$ for any $n$,
then $D^*(\Gamma)=D(\Gamma)= \{2\}$ with the exception that the groups listed in $(i)$
have $D^*(\Gamma)=\{2,3\}$, $D^*(L_3(2))=D^*(M_{11})=\{2,3,4\}$, and $D^*(M_{12})=\{2,4\}$.

$(i)$ $L_3(2).2$, $L_2(8)$, $L_2(8).3$, $L_2(9).22$, $L_2(9).23$, $L_2(9).22$, $L_2(11)$, $L_2(11).2$, $L_3(3)$, $L_2(13).2$, $L_2(16).2$, $L_2(16).4$, $L_3(4).3.22$, $M_{22}$, $M_{22}.2$, $M_{23}$, $M_{24}$.
\end{theorem}

Furthermore, we need the following preliminaries: 
Let $\Gamma$ be a group and $H\leqslant \Gamma$. For $\sigma \in \Gamma$, the set $\sigma H$ is called a \textit{left coset}  of $H$ in $\Gamma$. Similarly,
the set $H\sigma$ is called a \textit{right coset}  of $H$ in $\Gamma$. The set of all left cosets of $H$ in $\Gamma$ and the set of all right cosets have an equal cardinality. The set of distinct left (or right) cosets of $H$ in $\Gamma$ forms a partition of $\Gamma$; i.e., there exist $\sigma_1,\ldots , \sigma_k\in \Gamma \setminus H$ such that $\Gamma =H\cup \sigma_1H\cup \cdots \cup  \sigma_kH$ (equivalently, $\Gamma =H\cup H\sigma_1\cup \cdots \cup  H\sigma_k$). The index of $H\leqslant \Gamma$, denoted by $[\Gamma , H]$, is the cardinality of the set  of all distinct left (or right) cosets of $H$ in $\Gamma$. A subgroup $H$ of $\Gamma$ is normal in $\Gamma$, written $H\unlhd \Gamma$, if
$\sigma H=H\sigma$ for all $\sigma \in \Gamma$. If $H\unlhd \Gamma$ and $H\neq \Gamma$, we write $H\lhd \Gamma$.

\medskip

In Section 2, we obtain upper bounds for the distinguishing number of a set $X$ under the action of a group $\Gamma$ based on the distinguishing number of $X$ under the action of $H$, where $H$ is a subgroup of $\Gamma$. Also, we present some groups $\Gamma$ and subgroups $H$ of $\Gamma$ such that 
$\overline{D}(H)>\overline{D}(\Gamma)$, where $\overline{D}(\Gamma)={\rm max}\{D_{\Gamma}(X)~|~\Gamma ~{\rm acts~faithfully~ on}~X\}$. These examples give negative answer to a question of  Chan in \cite{Chan2}. Also we propose some cases for which the inequality $\overline{D}(H)\leq \overline{D}(\Gamma)$ is true.  
Motivated by a  problem 
(Chan 2006),  we characterize $D_{\Gamma,H}(X)$ which is the smallest number of labels admitting a labeling of $X$ such that the only elements of $\Gamma$ that induce label-preserving permutations lie in $H$, in Section 3.

\section{Upper bounds of $D_{\Gamma}(X)$ based on $D_{H}(X)$ }
We begin this section with the following question which  Potanka has asked   in \cite{Potanka}: 

\begin{question}\label{quesPot}
Let $\Gamma$ be a group and $H\leqslant \Gamma$. If  the distinguishing number of $X$ under the action of $H$ is $k$, then what can be said about the distinguishing
number of $X$ under the action of the entire $\Gamma$? 
\end{question} 

Note that if a given faithful action of $H$ on $X$ can be extended to a faithful action of $\Gamma$ on $X$ then $D_H(X) \leq D_{\Gamma}(X)$, because any labeling of $X$ that is distinguishing with respect to the action of $\Gamma$, is also distinguishing with respect to the action of $H$. However, is not true that,  every faithful action of $H$ on $X$ can necessarily be extended to a faithful action of $\Gamma$ (for example if $|\Gamma| > |X|!$). 

\medskip
 In this section,  we first answer to Question \ref{quesPot} when the action of any subgroup $H$ of $\Gamma$ on $X$ is the restriction of  action $\Gamma$ on $X$ to $H$, and next we state some results for that faithful action of $H$ on $X$ which cannot extended to a faithful action of $\Gamma$.

 Let $H$ be a proper  subgroup of $\Gamma$  such that  the action of any subgroup  of $\Gamma$ on $X$ is the restriction of  action $\Gamma$ on $X$ to that subgroup, so we can see that $2\leq D_H(X)\leq D_{\Gamma}(X)$. By Theorem \ref{abelidihedral} if $\Gamma$ is  abelian, then $ D_H(X)= D_{\Gamma}(X)=\{2\}$.  Since $2\leqslant D_H(X)\leq  D_{\Gamma}(X)$, so if $D_{\Gamma}(X) \leq 3$, then $D_{\Gamma}(X)\leq D_{H}(X)+1$. 

Now we obtain an upper bounds for $D_{\Gamma}(X)$ using the value $D_{H}(X)$ by the following theorem.
\begin{theorem}\label{thm}
Let $H_1,\ldots , H_k$ be all maximal nonidentity subgroups of $\Gamma$. If the action of any subgroup of $\Gamma$  on  $X$ is the restriction of the action of $\Gamma$ on $X$ to that subgroup, and   ${\rm max}\{D_{H_i}(X)\}_{i=1}^k=c$, then $D_{\Gamma}(X)\leq c+1$.
\end{theorem}
\proof  The proof is similar to the proof of Lemma 2.1 in \cite{Chan2}. Choose  a set $U$ of representatives of the orbits of $\Gamma$ on $X$, and let $L=\{\sigma\in \Gamma :~ u^{\sigma}=u ~{\rm for~each}~ u\in U\}$ stabilizes the set $U$ pointwise. We know that $X \setminus U$ is nonempty, because $\Gamma$ is nontrivial. Considering the action of $L$ on $X \setminus U$, we have $D_L(X\setminus U)\leq D_L(X)\leq c$, since $L$ is a subgroup of one of the $H_1,\ldots , H_k$. Suppose $C:X\setminus U \rightarrow \{1,2,\ldots , c\}$ is a $c$-labeling of $X\setminus U$  that is distinguishing with respect to the action of $L$. Now we define $C':X \rightarrow \{1,2,\ldots , c+1\}$ as 
\begin{equation*}
C'(x)=\left\{
\begin{array}{ll}
c+1 & {\rm if}~x\in U\\
C(x) & {\rm if}~x\notin U.
\end{array}\right.
\end{equation*}

We claim that $C'$ is a $(c+1)$-distinguishing labeling of $X$ with respect to the action of $\Gamma$. If $\sigma \in \Gamma$ and $\sigma$ preserves $C'$, then $\sigma$ should  fix each orbit representative $u\in U$, because they are only elements of label $c+1$ and are in different orbits, so $\sigma \in L$. Hence consider the action of $\sigma$ on $X\setminus U$. Since the restriction of $C'$ to $X\setminus U$  is a distinguishing labeling with respect to the action of $L$, and $\sigma \in L$ preserves this labeling, we have $\sigma$ is the identity. Therefore $C'$ is $(c+1)$-distinguishing labeling of $X$ with respect to the action of $\Gamma$.\qed

The bound presented in Theorem \ref{thm} is sharp. Let $M$ be the cycle graph $C_3$, and $\Gamma=S_3$ (the symmetric group on three letters). Each maximal subgroup of $\Gamma$ is abelian, and so the value of $c$ is two. Since $D_\Gamma(V(M))=3$, so this bound is sharp.
\begin{corollary}
Let  $\Gamma$ be a group acting on the set $X$. If  the action of any subgroup of $\Gamma$  on  $X$ is the restriction of the action of $\Gamma$ on $X$ to that subgroup, and $U$ contains the representatives of the orbits of $\Gamma$ on $X$, and $L=\{\sigma\in \Gamma :~ u^{\sigma}=u ~{\rm for~each}~ u\in U\}$ stabilizes the set $U$ pointwise, then $D_{\Gamma}(X)\leq D_L(X)+1$.
\end{corollary}
\proof Set $D_L(X)=c$. Suppose $C:X\setminus U \rightarrow \{1,2,\ldots , c\}$ is a $c$- labeling of $X\setminus U$  that is distinguishing with respect to the action of $L$. Now we define $C':X \rightarrow \{1,2,\ldots , c+1\}$ as 
\begin{equation*}
C'(x)=\left\{
\begin{array}{ll}
c+1 & {\rm if}~x\in U\\
C(x) & {\rm if}~x\notin U.
\end{array}\right.
\end{equation*}

We claim that $C'$ is a $(c+1)$-distinguishing labeling of $X$ with respect to the action of $\Gamma$. If $\sigma \in \Gamma$ and $\sigma$ preserves $C'$, then $\sigma$ should fix each orbit representative $u\in U$, because they are only elements of label $c+1$ and are in different orbits, so $\sigma\in L$. Hence consider the action of $\sigma$ on $X\setminus U$. Since the restriction of $C'$ to $X\setminus U$  is a distinguishing labeling with respect to the action of $L$, and $\sigma\in L$ preserves this labeling, so  $\sigma$ is the identity. Therefore $C'$ is $(c+1)$-distinguishing labeling of $X$ with respect to the action of $\Gamma$. \qed

\medskip
Now  by the concept of motion, we obtain an upper bound for $D_{\Gamma}(X)$ based on $D_{H}(X)$,  where  the faithful action of $H$ on $X$ cannot necessarily be extended to a faithful action of $\Gamma$.
Let $\Gamma$ be a finite group acting on the set $X$. For $\sigma \in \Gamma$ let $m(\sigma) = \{x \in X :~ \sigma(x) \neq x\}$ and let $m(X) = {\rm min}\{m(\sigma) :~ \sigma \neq id\}$. Call $m(\sigma)$ the motion of $\sigma$ and $m(X)$ the motion of $X$. We state the following definition: 

\begin{definition}
Let $H$ be a subgroup of  $\Gamma$. Suppose that $C$ is a distinguishing labeling of $X$ under $H$ with $D_H(X)$ labels, and  $\sigma_i$ ($1\leq i\leq k$) are all elements of $\Gamma$ which preserve the labeling $C$.   The partition   $T=\{I_1,\ldots , I_n\}$   of the set $\{1,\ldots , k\}$ is called a {\rm good partition}, if $\bigcap_{i\in I_j}m(\sigma_i)\neq \emptyset$ for every $j=1,\ldots ,n$.   
\end{definition}

\begin{theorem}\label{thm1}
Let $X$ be a set, and $\Gamma$ and $H$ be two groups acting on $X$ where $H\leqslant \Gamma$.   Let $t={\rm min}\big\{|T|:~T ~{\rm is~ a~ good~ partition~of}~\{1,\ldots , k\} \big\}$. We have 
\begin{enumerate}
	\item [(i)] $D_{\Gamma}(X)\leq D_H(X)+t$.
	\item[(ii)]  If $ \bigcap_{i=1}^k m(\sigma_i)\neq \emptyset$, then $D_{\Gamma}(X)\leq D_H(X)+1$.	
\end{enumerate} 
\end{theorem}
\proof 
\begin{enumerate}
	\item [(i)] Let $T$ be a good partition of $\{1,\ldots , k\}$  with $|T|=t$. Set $T=\{I_1,\ldots , I_t\}$. With respect to minimality of  $t$, we have
$$\left(   \bigcap_{i\in I_{k'}} m(\sigma_i)  \right) \cap \left(   \bigcap_{i\in I_{k''}} m(\sigma_i)  \right)=\emptyset,$$
for $k'\neq k''$. Let $C$ be a distinguishing labeling of $X$ under $H$ with $D_H(X)$ labels. We choose exactly one element from each $ \bigcap_{i\in I_j} m(\sigma_i)$ ($j=1,\ldots , t$), say $x_j$, and change its  label to the label $D_H(X)+j$.  It can be seen that this new labeling is not preserved by elements  $\sigma_1, \ldots , \sigma_k$, and so the new labeling is distinguishing under  $\Gamma$. Therefore we obtained a distinguishing labeling of $X$ under $\Gamma$ with $D_H(X)+t$ labels.

\item[(ii)] It follows from Part (i). \qed
\end{enumerate} 

\begin{corollary}
Let $X$ be a set, and $\Gamma$ and $H$ be two groups acting on $X$ where $H\leqslant \Gamma$. Suppose that $D_H(X)=d$, and $D(X,H)$ is the set of all nonisomorphic $d$-distinguishing labeling of $X$ under $H$. We denote the elements of $D(X,H)$ by $\phi_i$ where $1\leq i \leq |D(X,H)|$. Let $\sigma^{(i)}_1 , \ldots ,\sigma^{(i)}_{k_i}$ be the elements of $\Gamma$ which preserve the labeling $\phi_i$, and  $t_i={\rm min}\big\{|T|:~T ~{\rm is~ a~good~ partition~of}~\{1,\ldots , k_i\}\big\}$. If $t={\rm min}\{t_i\}_{i=1}^{|D(X,H)|}$, then $D_{\Gamma}(X)\leq D_H(X)+t$.
\end{corollary}

If $H=\{1\}$, then we can obtain an upper bound for the value of $D_{\Gamma}(X)$ by Theorem \ref{thm1} as follows: 
\begin{corollary}
If $\Gamma$ is a group of order $k+1$ acting on the set $X$, then $D_{\Gamma}(X)\leq 1+t$ where $t={\rm min}\big\{|T|:~T ~{\rm is~ a~ good~ partition~of}~\{1,\ldots ,k\}\big\}$.
\end{corollary}

We continue  this section with the following  problem from  Chan \cite{Chan2}.

\begin{problem} 
Let $\Gamma$ be a finite group, and $\overline{D}(\Gamma)={\rm max}\{D_{\Gamma}(X)~|~\Gamma ~{\rm acts~faithfully~ on}~X\}$. If $H$ is a subgroup of $\Gamma$, dose it follow that $\overline{D}(H)\leq \overline{D}(\Gamma)$? 
\end{problem}

In general, the answer of this problem is  negative. For instance, we use the dihedral groups. If $\Gamma=D_n$, and $H=D_d$ with $n\geq 12$ and $d\in \{3,4,5,6,10\}$ such that $d|n$, then $\overline{D}(\Gamma)=2$ and $\overline{D}(H)=3$, by Theorem \ref{abelidihedral}.  It is clear that if $\Gamma$ is a group and the symmetric group $S_n$ is its subgroup such that $\overline{D}(\Gamma)\leq n-1$, then by Theorem \ref{symmthm} $\overline{D}(\Gamma) < \overline{D}(S_n)=n$. As another  example, let $\Gamma=M_{11}$. It is known that $H=S_5$ is a maximal subgroup of $\Gamma$. Now by Theorems \ref{symmthm} and \ref{almostthm}   we have $4=\overline{D}(\Gamma) < \overline{D}(H)=5$. By a similar argument, if $\Gamma$ is a group and the alternative group $A_n$ is its subgroup such that $\overline{D}(\Gamma)\leq n-2$, then $\overline{D}(\Gamma) < \overline{D}(A_n)=n-1$.

\medskip

Note that  there exist cases for which $\overline{D}(H)\leq \overline{D}(\Gamma)$ is true. In the rest of this section, we characterize cases for which $\overline{D}(H)\leq \overline{D}(\Gamma)$. First we state the following result which is an immediate consequence of Theorem \ref{abelidihedral}:
\begin{theorem}
If $\Gamma$ is an abelian group, and $H$ is a subgroup of $\Gamma$, then $\overline{D}(H)\leq \overline{D}(\Gamma)$.
\end{theorem}

\begin{theorem}
 Let $\overline{D}(H) =c< \infty$,   and $A$ be the set of all pairs $(X,H)$ for which $D_H(X)=c$. If at least one of the pairs in $A$, say $(X,H)$, can be extended to a faithful action of $\Gamma$, then $\overline{D}(H)\leq \overline{D}(\Gamma)$. 
\end{theorem}
\proof
Since the action of the pair $(X,H)$ from $A$   can be extended to a faithful action of $\Gamma$, so $D_H(X)\leq D_{\Gamma}(X)$. Now   using $D_{\Gamma}(X) \leq \overline{D}(\Gamma)$ and $D_H(X)=\overline{D}(H)$, we  conclude that  $\overline{D}(H)\leq \overline{D}(\Gamma)$.\qed

The following theorem is a direct consequence of Theorem \ref{uppdisgroup}.
\begin{theorem}\label{lessthm}
Let $\Gamma$ and $H$ be two groups and $H$ be a subgroup of $\Gamma$. 
\begin{itemize}
\item[(i)] If $|H|\leq (\overline{D}(\Gamma))!$, then $\overline{D}(H)\leq \overline{D}(\Gamma)$.
\item[(ii)] If $|\Gamma|\leq (\overline{D}(H))!$, then $\overline{D}(\Gamma)\leq \overline{D}(H)$.
\end{itemize}
\end{theorem}

\begin{corollary}
If $\Gamma=S_n$ or $\Gamma=A_n$ where $n\geq 3$, then for every subgroup $H$ of $\Gamma$, $\overline{D}(H)\leq \overline{D}(\Gamma)$.
\end{corollary}
\proof   If follows from Theorems \ref{symmthm} and \ref{lessthm}.\qed

\begin{theorem}
Let  $\Gamma$ be a group such that $|\Gamma|\leq 2(\overline{D}(\Gamma)!)$. If $H$ is a subgroup of  $\Gamma$, then $\overline{D}(H)\leq \overline{D}(\Gamma)$.
\end{theorem}
\proof If $H$ is a trivial group or $H=\Gamma$, then the result is clear. So we suppose that $H$ is a nontrivial proper subgroup of  $\Gamma$. By  contradiction suppose  that  $\overline{D}(H) > \overline{D}(\Gamma)$. Hence $|H| > \overline{D}(\Gamma)!$, by Theorem \ref{lessthm}. Since $[\Gamma , H]\geq 2$, we can conclude that $|\Gamma| > 2 \overline{D}(\Gamma)!$,  which is a contradiction.\qed

This result can easily be generalized   as follows:
\begin{corollary}
Let  $\Gamma$ be a group such that for an integer number $t$, $|\Gamma|\leq t(\overline{D}(\Gamma)!)$. If $H$ is a subgroup of  $\Gamma$ such that $[\Gamma , H]\geq t$, then $\overline{D}(H)\leq \overline{D}(\Gamma)$.
\end{corollary}

A theorem by Gluck \cite{Gluck} states $D_{\Gamma}(X)=2$ whenever $|\Gamma|$ is odd, and hence $\overline{D}(\Gamma)=2$. Our hypothesis on $|\Gamma|$
implies  the two following results directly.
\begin{theorem}
\begin{itemize}
\item[(i)] Let $\Gamma$ be a group of odd order. If $H$ is a subgroup of $\Gamma$, then $\overline{D}(H)\leq \overline{D}(\Gamma)$.
\item[(ii)] Let $\Gamma$ be a nontrivial group. If $H$ is a subgroup of odd order of $\Gamma$, then $\overline{D}(H)\leq \overline{D}(\Gamma)$.
\end{itemize}
\end{theorem}

\section{Characterization of $D_{\Gamma ,H}(X)$}

Chan in \cite{Chan2} used subgroup structure as a method  to generalize the notion of the
distinguishing number, as follows. Given a group $\Gamma$ acting faithfully on a set $X$ and $H$ a subgroup of $\Gamma$, let $D_{\Gamma ,H}(X)$ denote the smallest number of labels admitting a labeling of $X$ such that the only elements of $\Gamma$ that induce label-preserving permutations lie in $H$. Thus, when $H = 1$, we recover the original notion of the distinguishing number. It is clear that $D_{\Gamma, H}(X)\leq D_{\Gamma}(X)$, $D_{\Gamma, 1}(X)= D_{\Gamma}(X)$, and $D_{\Gamma, \Gamma}(X)=1$. He proposed the following problem in \cite{Chan2}.

\begin{problem} 
 Characterize $D_{\Gamma ,H}(X)$.
 \end{problem}

In this section we shall characterize $D_{\Gamma ,H}(X)$. To do this, we begin with the following lemma:  
\begin{lemma}\label{lemm1}
Let $\Gamma$ be a group,  $H_1$ and $H_2$ be two subgroups of $\Gamma$. If $H_1\leqslant H_2$, then $D_{\Gamma, H_2}(X)\leq D_{\Gamma , H_1}(X)$.
\end{lemma}
\proof  Set $D_{\Gamma , H_1}(X)=d$. Suppose that $\phi: X\rightarrow \{1,\ldots , d\}$ is a $d$-labeling of $X$, such that the only elements of $\Gamma$ that induce label-preserving permutation lie in $H_1$. Since $H_1\leqslant H_2$, so each element of $\Gamma$ that induce label-preserving permutation lie in $H_2$. Therefore $D_{\Gamma, H_2}(X)\leq d$.\qed

\begin{corollary}
If $\Gamma$ is a group, and $\{e\}\leqslant H_1 \leqslant H_2\leqslant \cdots \leqslant H_k \leqslant \Gamma$ is  a chain of subgroups of $\Gamma$, then 
\begin{equation*}
1\leq D_{\Gamma, H_k}(X)\leq D_{\Gamma , H_{k-1}}(X)\leq \cdots \leq D_{\Gamma, H_1}(X)\leq D_{\Gamma}(X).
\end{equation*}
\end{corollary}

\begin{lemma}\label{lemm10}
Let $\Gamma_1$ and $\Gamma_2$ be  two groups such that $\Gamma_1\leqslant \Gamma_2$, and $H\leqslant \Gamma_1$, then $D_{\Gamma_1, H}(X)\leq D_{\Gamma_2 , H}(X)$.
\end{lemma}
\proof  Set $D_{\Gamma_2 , H}(X)=d$. We suppose that $\phi: X\rightarrow \{1,\ldots , d\}$ is a $d$-labeling of $X$, such that the only elements of $\Gamma_2$ that induce label-preserving permutation lie in $H$. Since $\Gamma_1\leqslant \Gamma_2$, so each element of $\Gamma_1$ that induce label-preserving permutation lie in $H$. Therefore $D_{\Gamma_1, H}(X)\leq d$.\qed

\begin{theorem}\label{normalthm}
If $\Gamma$ is a group and $H$ is its normal subgroup, $H\unlhd \Gamma$ (the action of $\sigma H\in \Gamma /H$ on $X$ is as $\sigma H (x)=\sigma (x)$, for all $x\in X$), then $D_{\Gamma, H}(X)= D_{\Gamma /H}(X)$. 
\end{theorem}
\proof  Set $D_{\Gamma , H}(X)=d$, and  $\phi: X\rightarrow \{1,\ldots , d\}$ be a $d$-labeling of $X$, such that the only elements of $\Gamma$ that induce label-preserving permutation lie in $H$. We want to show that $\phi$ is a $d$-distinguishing labeling of $X$ under $\Gamma /H$. If $\sigma H\in \Gamma /H- \{H\}$ is an element of $\Gamma /H$ preserving $\phi$, then $\phi (\sigma H (x))= \phi (x)$ for any $x\in X$. Thus $\phi (\sigma  (x))= \phi (x)$ for any $x\in X$. Since  the only elements of $\Gamma$ that induce label-preserving permutation lie in $H$, so $\sigma \in H$, which is a contradiction. Therefore $D_{\Gamma /H}(X)\leq D_{\Gamma, H}(X)$. For converse, let set $D_{\Gamma / H}(X)=t$, and  $\phi: X\rightarrow \{1,\ldots , t\}$ be a $t$-distinguishing labeling of $X$ under $\Gamma /H$. We want to show that $\phi$ is a $d$-labeling of $X$ such that the only elements of $\Gamma$  induce label-preserving permutation lie in $H$. 
 If $\sigma \in \Gamma$ is an element of $\Gamma$ preserving $\phi$, then $\phi (\sigma (x))= \phi (x)$ for any $x\in X$. Thus $\phi (\sigma H (x))= \phi (x)$ for any $x\in X$, and so $\sigma H = H$, which means $\sigma \in H$. Therefore $D_{\Gamma , H}(X)\leq D_{\Gamma /H}(X)$. \qed

\begin{corollary}
If  $\Gamma$ is a group and $H$ is its  subgroup  with $[\Gamma:H]=p$ such that $p$ is the least prime number with $p~| ~|\Gamma|$, then $D_{\Gamma , H}(X)=2$.
\end{corollary}
\proof By hypotheses we conclude that $H$ is a normal subgroup of $\Gamma$, and $\Gamma /H$ is an abelian group. So the result follows by Theorems \ref{abelidihedral} and \ref{normalthm}.\qed

\begin{corollary}
If  $\Gamma$ is a group and $H$ is a normal subgroup of $\Gamma$ such that $\Gamma /H$ is a nontrivial abelian group, then  $D_{\Gamma , H}(X)=2$, especially, $D_{\Gamma , \Gamma'}(X)=2$ where $\Gamma'$ is the derived subgroup of $\Gamma$.
\end{corollary}

\begin{corollary}
Let $\Gamma$ be a group and $H$ be a subgroup of $\Gamma$. If $<H>_{nor}=\bigcap_{H\leqslant K\unlhd \Gamma}K$, then $D_{\Gamma/<H>_{nor}}(X)\leq D_{\Gamma ,H}(X)$.
\end{corollary}
\proof  It can be shown that $H \leqslant <H>_{nor}\unlhd \Gamma$. So by Lemma \ref{lemm10}, $D_{\Gamma, H}(X)\leq D_{\Gamma ,<H>_{nor}}(X)$. Now we have the result by Theorem \ref{normalthm}.\qed

A group $\Gamma$ is called {\it metacyclic} if it has a normal subgroup $N \lhd \Gamma$ such that both $N$ and $\Gamma /N$ are cyclic. Such groups have been completely classified in [9], and include all groups of square free order.
\begin{corollary}{\rm \cite{Chan2}}
 Let $\Gamma$ be a metacyclic group. Then ${\rm max}  D^*(\Gamma) \leq 3$.
 \end{corollary}
We obtain as a special case that if $\Gamma$ is dihedral then ${\rm max}  D^*(\Gamma) \leq 3$, as shown in \cite{Albert}.

\begin{corollary}
\begin{itemize}
\item[(i)] If $H$ is a nontrivial proper subgroup of an abelian group $\Gamma$, then $D_{\Gamma , H}(X)=2$.
\item[(ii)] If $H$ is a subgroup of metacyclic group $\Gamma$, then $D_{\Gamma , H}(X) \leq 3$.
\end{itemize}
\end{corollary}

Let $\Gamma$ be a group acting faithfully on the set $X$, and $\phi$ be a  $d$-labeling of $X$. Suppose that $(\Gamma,\phi)$ is the set of all elements of $\Gamma$ preserving the labeling $\phi$. It can be seen that $(\Gamma,\phi)$ is a subgroup of $\Gamma$, and $D_{\Gamma, (\Gamma,\phi)}(X)=d$. If $H$ is a subgroup of $\Gamma$, and $D_{\Gamma ,H}(X)=c$, then $(\Gamma,\phi)\leqslant H$ where $\phi : X\rightarrow \{1,\ldots , c\}$ is a $c$-labeling of $X$. Hence
\begin{align*}
D_{\Gamma , H}(X)\leq {\rm min}\{d~|~\phi~{\rm is~a~}d-{\rm labeling~of}~ X~{\rm such~ that}~(\Gamma, \phi)\leqslant H\},\\
D_{\Gamma , H}(X)\geq {\rm max}\{d~|~\phi~{\rm is~a~}d-{\rm labeling~of}~ X~{\rm such~ that}~H\leqslant (\Gamma, \phi)\}.
\end{align*}

By this argument we can state the two following algorithms for obtaining an upper and a lower bound for $D_{\Gamma , H}(X)$. 
  Here, $L_d$ is the set of all $d$-labeling of $X$. Since $|L_d|= \sum_{j=0}^{|X|}(-1)^j {d \choose j} (d-j)^{|X|}$ (the number of  function from $X$ onto $\{1,\ldots ,d\}$), so we denote the elements of $L_d$ by $\phi_d^i$ where $1\leq i \leq |L_d|$.\\

Here we state an algorithm for obtaining an upper bound for $D_{\Gamma , H}(X)$:
\vspace{.6cm}

\begin{tabular}{l}

Input: $X,\Gamma, H$\\
Output: An upper bound for $D_{\Gamma , H}(X)$\\
Set $\mathfrak{D}\leftarrow \emptyset$\\
for($d=1, d\leq |X|, d++$)\\
\{\\
\indent Construct $L_d$\\
\indent for($i=1, i\leq |L_d|, i++$)\\
\indent \indent \indent \{\\
\indent \indent \indent If ($(G, \phi_d^i)\leqslant H$)\\
\indent \indent \indent Set $\mathfrak{D}\leftarrow \mathfrak{D}\cup d$\\
\indent \indent \indent \}\\
\}\\
Print $D_{\Gamma , H}(X)\leq {\rm min} \mathfrak{D}$ 
\end{tabular}

\vspace{.7cm}

Now we propose an algorithm for obtaining a lower bound for $D_{\Gamma , H}(X)$: 
\vspace{.4cm}

\begin{tabular}{l}

Input: $X,\Gamma, H$\\
Output: A lower bound for $D_{\Gamma , H}(X)$\\
Set $\mathfrak{D}\leftarrow \emptyset$\\
for($d=|X|, d\geq 1, d--$)\\
\{\\
\indent Construct $L_d$\\
\indent for($i=1, i\leq |L_d|, i++$)\\
\indent \indent \indent \{\\
\indent \indent \indent If ($H \leqslant (G, \phi_d^i)$)\\
\indent \indent \indent Set $\mathfrak{D}\leftarrow \mathfrak{D}\cup d$\\
\indent \indent \indent \}\\
\}\\
Print $D_{\Gamma , H}(X)\geq {\rm max} \mathfrak{D}$

\end{tabular}



\begin{corollary}
If $H$ and $K$ are two normal subgroups of $HK$ such that $H\cap K=\{1\}$, then $D_{HK,H}(X)= D_{K}(X)$.
\end{corollary}
\proof Since $H\cap K=\{1\}$, so  $HK/H\cong K$. Now by Lemma \ref{lemm10} we have the result.\qed

\begin{proposition}
If $H$ and $K$ are two subgroups of the group $HK$ (i.e., $HK=KH$), then $D_{HK,H}(X) \geq D_{K, H\cap K}(X)$.
\end{proposition}
\proof  Set $D_{HK , H}(X)=c$, and  $\phi: X\rightarrow \{1,\ldots , c\}$ be a $c$-labeling of $X$, such that the only elements of $HK$ that induce label-preserving permutation lie in $H$. We want to show that $\phi$ is a $d$-distinguishing labeling of $X$ under $\Gamma /H$. If $k\in K$  preserves the labeling  $\phi$, then $k\in H$, and so $k\in H\cap K$. Thus $\phi$ is a $c$-labeling of $X$, such that the only elements of $K$ that induce label-preserving permutation lie in $H\cap K$. Therefore $D_{HK , H}(X)\geq D_{K, H\cap K}(X)$.\qed

\begin{corollary}
If $H$ and $K$ are two groups of the group $HK$ such that $H\cap K=\{1\}$, then $D_{HK,H}(X) \geq D_{K}(X)$.
\end{corollary}

\section{Conclusion}

In this paper we gave an upper bound on the distinguishing number of a set knowing the distinguishing number of a set under the action of a subgroup.  Also, we gave some groups $\Gamma$ and subgroups $H$ of $\Gamma$ such that 
$\overline{D}(H)>\overline{D}(\Gamma)$.  We proposed some cases for which the inequality $\overline{D}(H)\leq \overline{D}(\Gamma)$ is true.  
Also motivated by a  problem 
(Chan 2006),  we characterized $D_{\Gamma,H}(X)$.  
 Finally, we stated  two algorithms for obtaining an upper and a lower bound for $D_{\Gamma , H}(X)$. We end the paper with proposing the following question:  

\begin{question}
Let $\Gamma$ be a group acting on the sets $X$ and $Y$, where $Y  \subseteq X$. If  the distinguishing number of $Y$ under the action of $\Gamma$ is $k$, then what can be said about the distinguishing number of $X$ under the action of $\Gamma$? 
\end{question} 


\end{document}